\documentclass[MSNbibl,number,citesort,seceqn,dvips]{arxbj}
\usepackage{upgreek}
\usepackage{mathbh}

\aid{0}
\volume{19}
\issue{2}
\pubyear{2013}
\firstpage{363}
\lastpage{386}
\doi{10.3150/11-BEJ402}

\makeatletter
\newtheorem{Lemma}{Lemma}[section]
\newtheorem{Proposition}[Lemma]{Proposition}
\newtheorem{Theorem}[Lemma]{Theorem}

\newtheorem{Corollary}[Lemma]{Corollary}
\newcommand{\IIst}{\mathcal{I}^*}
\newcommand{\prob}{\mathbb{P}}
\newcommand{\PP}{\mathcal{P}}
\newcommand{\GG}{\mathcal{G}}
\newcommand{\Rbold}{{\mathbb{R}}}
\newcommand{\Zbold}{{\mathbb{Z}}}

\newcommand{\expec}{\mathbb{E}}
\newcommand{\SSS}{\mathcal{S}}
\newcommand{\SSst}{\mathcal{S}^*}
\newcommand{\Op}{\mathrm{O}_{ \prob}}
\newcommand{\convd}{\stackrel{d}{\longrightarrow}}
\newcommand{\convp}{\stackrel{ {\mathbb P}}{\longrightarrow}}
\newcommand{\vep}{\varepsilon}
\newcommand{\e}{\mathrm{e}}
\newcommand{\EXP}{\operatorname{Exp}}
\newcommand{\Poi}{\operatorname{Poi}}
\newcommand{\eqref}[1]{(\ref{#1})}
\makeatother

\begin{document}
\begin{frontmatter}

\title{Weak disorder in the stochastic mean-field model of distance II}
\runtitle{Distances in weak disorder}

\begin{aug}
\author[1]{\fnms{Shankar} \snm{Bhamidi}\thanksref{1}\ead[label=e1]{bhamidi@email.unc.edu}},
\author[2]{\fnms{Remco}~\snm{van~der~Hofstad}\thanksref{2}\ead[label=e2]{rhofstad@win.tue.nl}}
\and\break
\author[3]{\fnms{Gerard}~\snm{Hooghiemstra}\corref{}\thanksref{3}\ead[label=e3]{g.hooghiemstra@tudelft.nl}}
\runauthor{S. Bhamidi, R. van der Hofstad and G. Hooghiemstra} 
\address[1]{Department of Statistics and Operations Research,
The University of North Carolina, 304 Hanes Hall, Chapel Hill, NC
27510, USA. \printead{e1}}
\address[2]{Department of Mathematics and
Computer Science, Eindhoven University of Technology, P.O.~Box 513,
5600 MB Eindhoven, The Netherlands. \printead{e2}}
\address[3]{DIAM, Delft University of Technology, Mekelweg 4, 2628CD
Delft, The
Netherlands.\\
\printead{e3}}
\end{aug}

\received{\smonth{9} \syear{2010}}
\revised{\smonth{8} \syear{2011}}

%
\begin{abstract}
In this paper, we study the complete graph $K_n$ with $n$ vertices,
where we attach
an independent and identically distributed (i.i.d.) weight to each of
the $n(n-1)/2$ edges. We focus on the weight $W_n$ and the number of
edges $H_n$ of the minimal weight path between
vertex $1$ and vertex $n$.

It is shown in (\textit{Ann. Appl. Probab.} \textbf{22} (2012) 29--69) that when
the weights on the edges are
i.i.d. with distribution
equal to that of $E^s$, where $s>0$ is some parameter, and $E$ has an
exponential distribution with mean 1, then $H_n$ is asymptotically
normal with asymptotic mean $s\log n$
and asymptotic variance $s^2\log n$. In this paper,
we analyze the situation when the weights have distribution $E^{-s},
s>0$,
in which case the behavior of $H_n$ is markedly different as $H_n$ is a
tight sequence of random variables. More precisely, we
use the method of Stein--Chen for Poisson approximations to show
that, for almost all $s>0$, the hopcount $H_n$ converges
in probability to the nearest integer of $s+1$
greater than or equal to 2, and identify the limiting
distribution of the recentered and rescaled minimal weight.
For a countable set of special $s$ values denoted by
$\mathcal{S}=\{s_j\}_{j\geq2}$, the hopcount $H_n$ takes on the values
$j$ and $j+1$ each with \textit{positive} probability.
\end{abstract}

%

\end{frontmatter}

\section{Introduction}
\label{sec-intro}
One of the central themes of modern discrete probability is
the study of the effect of random edge
disorder on various properties of the underlying network.
The base network itself could be \textit{deterministic}, for example, a
large finite box in the lattice or the
complete graph on $n$ vertices, or \textit{random}, for example,
the giant component of the Erd\H{o}s--R\'enyi
random graph or the configuration model. Each edge is assigned a random
edge weight, whose interpretation varies depending on the context. One
can think of this weight as the \textit{cost} in traversing the edge,
yielding first passage percolation-type models.
Alternatively, one could think of the underlying graph as an electrical
network and the assigned weights as resistances, yielding a random
resistor network, or as capacities on edges and the underlying graph as
a flow carrying network, entrusted with carrying flow
(commodities, information, etc.) between various parts of the network.

One graph model that has resulted in many problems of fundamental
interest is the complete graph $K_n$ on $n$
vertices with random edge weights. In the various contexts mentioned
above, this model both gives rise to
very interesting conjectures as well as generates new techniques and
insights in
probability theory that can then be applied in a number of other contexts.
While providing a complete list of references of the various models
that have been studied in this context
would be impractical, we direct the interested reader to \cite
{janson1999one} for one of the first refined
results in first passage percolation in this context, \cite
{frieze1985value} for an analysis of
the cost of the minimal spanning tree, \cite{grimmett1984random}
for an analysis of random electrical networks on the complete graph,
\cite{aldous20011} for a study of the random assignment problem, \cite
{aldous2009uniform} for an
analysis of the multicommodity flow problem
and the survey paper \cite{aldous2003objective}, where a number of
other examples are analyzed via the powerful \textit{local weak
convergence} method.

Let us now focus on the particular problem dealt with in this study and
the motivations behind it.
Suppose we start with a connected graph $\GG_n$ (deterministic such as
$K_n$ or random) on $n$ vertices.
Suppose each edge $e$ is assigned a random positive edge weight $E_e$.
We shall assume that the weights are i.i.d.~over the edges with
some distribution $F$, with density $f$.
Fix two vertices (say chosen uniformly at random from $\GG_n$), and let
us denote them by $1$ and $n$.
For any path $\PP$ between the two vertices, let the weight of the path
$w(\PP)$ be defined by
\[
w(\PP):= \sum_{e\in\PP} E_e,
\]
that is, the sum of weights of the edges in the path. The
\textit{optimal} or \textit{minimal weight} path (which is unique
since the edge
weights have a density) is the path that minimizes the above weight
function. In the study of random systems, this regime is often called
the \textit{weak disorder} regime, while probabilists know this
problem as
``first passage percolation.'' The mental picture one can have is that
the network is entrusted with carrying flow between various
nodes of the network, and the way it performs this duty is via routing
flow through optimal paths. We shall defer a more extensive discussion
of the relevant literature to Section~\ref{sec:disc}.

Another regime which is of tremendous interest is the \textit{strong
disorder regime}.
Here the weight of a path is either the \textit{maximum} or the
\textit{minimum} weight of all edges in the path. We denote the weight
functions as
%
%
\begin{equation}\label{eqn:wmax-model}
w_{\mathtt{max}}(\PP):= \max_{e\in\PP} E_e,
\end{equation}
and
%
%
\begin{equation}
\label{eqn:wmin-model}
w_{\mathtt{min}}(\PP):= \min_{e\in\PP} E_e.
\end{equation}
In both situations, one is interested in properties of the path which
minimizes the above weight
function. One is also interested in formulating a model, depending on a
real-valued parameter,
the ``inverse parameter,'' which interpolates between these two models.
One can then study questions such as phase transitions, where there is
a change in the behavior
of the system from the weak disorder regime to the strong disorder regime.
Given the set of edge weights $E_e$, one method of doing this is as
follows: assign each edge a cost
$E_e^{\beta}$ where $\beta\in\Rbold$ is a real-valued parameter.
With these edge weights, suppose that, as before,
we consider the weak disorder regime, so that now the weight of a path
$\PP$ is
\[
w_{\beta}(\PP):= \sum_{e\in\PP} E_e^{\beta}.
\]
Then, we can identify the following special cases:
\begin{longlist}[(a)]
\item[(a)] \textit{Original model}: $\beta=1$ is our original model.
\item[(b)] \textit{Graph distance}: $\beta=0$ gives us the graph distance
between the chosen vertices in the graph $\GG_n$.
\item[(c)] \textit{Strong disorder, max edge weight}: The case $\beta\to
+\infty
$ gives us the strong disorder regime where the weight of a path is
given by \eqref{eqn:wmax-model}. This is also called the \textit{minimal
spanning tree} regime as the optimal path between the two vertices is
the same as the path in the minimal spanning tree on $\GG_n$ with edge
weights $E_e$.
\item[(d)] \textit{Strong disorder, min edge weight}: $\beta\to-\infty$ gives
us the strong disorder model where the weight of a path is given by
\eqref{eqn:wmin-model}.
\end{longlist}

Thus this model allows us to interpolate between various regimes of
interest. We shall denote the optimal path by $\PP_{\mathtt{opt}}(\beta)$.
Given a particular base network $\GG_n$ and edge weight distribution~$E_e$, two statistics are of paramount interest:
\begin{longlist}[(ii)]
\item[(i)]\textit{Minimal weight:} This is the actual weight of the optimal
path, namely $W_n =  \sum_{e\in\PP_{\mathtt{opt}}(\beta)} E_e^{\beta}$.\vspace*{1.5pt}
\item[(ii)]\textit{Hopcount:} This is defined as the number of edges in the
optimal path $\PP_{\mathtt{opt}}(\beta)$. We shall denote this random
variable by $H_n(\beta)$.
\end{longlist}

\textit{Aim of this paper}:
In this paper, we shall specialize to the case where the graph $\GG_n$
is the complete graph $K_n$ and
each edge originally has edge weight $E_e^{\beta}$, where
$E_e$ is exponentially distributed with rate $1$ ($E_e \stackrel{d}{=}
\EXP(1)$).
We shall study the case where $\beta< 0$. The case where $\beta>0$ has
been solved in \cite{BH09}, where it was proved that, for $\beta> 0$,
%
%
\begin{equation}
\label{eqn:hns-s-big-0}
\frac{H_n(\beta)-\beta\log{n}}{\sqrt{\beta^2\log{n}}} \convd Z,
\end{equation}
where $Z$ denotes a standard normal random variable, and $\convd$ denotes
convergence in distribution.
In the same paper it was proved that, for the optimal weight
$W_n=W_n(\beta)$, there exists a constant $\lambda=\lambda(\beta)> 0$
and a non-degenerate real-valued random variable $\Xi(\beta)$ such that
%
%
\begin{equation}\label{eqn:wn-s-big-0}
W_n(\beta)-\frac{1}{\lambda}\log{n} \convd\Xi(\beta).
\end{equation}
In this study, we shall derive asymptotics for the two random variables
of interest $W_n(\beta)$ and $H_n(\beta)$ as $n\to\infty$, and see
that the behavior in the case when $\beta<0$ is markedly different.

Throughout the paper, we make use of the following standard notation.
We let $\convd$ denote convergence in distribution, and
$\convp$ convergence in probability. For a sequence of random variables
$(X_n)_{n\geq1}$, we write $X_n=\Op(1)$ when $|X_n|$
is a tight sequence of random variables as $n\rightarrow\infty$, and
$X_n=\mathrm{o}_{ \prob}(1)$ when $|X_n|\convp0$ as $n\rightarrow\infty$.
For a non-negative function $n\mapsto g(n)$,
we write $f(n)=\mathrm{O}(g(n))$ when $|f(n)|/g(n)$ is uniformly bounded, and
$f(n)=\mathrm{o}(g(n))$ when $\lim_{n\rightarrow\infty} f(n)/g(n)=0$.
We let $\EXP(\lambda)$ denote an exponential random variable with
rate $\lambda$ and $\Poi(\lambda)$ a Poisson random variable with
mean $\lambda$. We write that a sequence of events $(\mathcal{E}_n)_{n\geq1}$
occurs \textit{with high probability} (\textit{w.h.p.}) when
$\prob(\mathcal{E}_n)\rightarrow1$.
Finally, for $x\in{\mathbb R}$,
we denote by $\lfloor x\rfloor$ the largest integer smaller than or
equal to $x$ and
by $\lceil x\rceil$ the smallest integer larger than or equal to~$x$.

We now state our main results and defer a further discussion to
Section~\ref{sec:disc}.

\section{Results}
\label{sec:results}
Before stating the main result, we need some further notation.
We study the complete graph~$K_n$ with i.i.d. edge weights $E^{-s}_{\{
i,j\}}, 1\le i< j\le n$,
on the edges of $K_n$. Thus, compared to the discussion
in the previous section, we have taken $s=-\beta$, and we shall
study the $s>0$ regime. For fixed $s>0$, define the function
%
\begin{equation} \label{def_g} g_s(x)=\frac{x^{s+1}}{(x-1)^s},\qquad x\ge
2.
\end{equation}
Observe that, for $0<s\le1$, the function $g_s(x), x\ge2,$ is increasing,
while for $s>1$, the function
is strictly convex with unique minimum at $x=s+1$.
We shall be interested in minimizing this function only on the set
$\Zbold_+$ of positive integers.
Then there is a sequence of
values $s=s_j, j\ge2,$ for which the minimum integer of $g_{s}$ is
\textit{not} unique.
From the equation $g_s(j)=g_s(j+1)$, and the bounds $j-1<s<j$, it is
not hard to verify that
%
\begin{equation} \label{seq_sj} s_j=\frac{\log(1+j^{-1})}{\log
(1+(j^2-1)^{-1})}\in(j-1,j),\qquad j=2,3,\ldots.
\end{equation}
We will need to deal with these special points separately. When
$s\notin\SSS=\{s_2,s_3,\ldots\}$, then there is a unique
integer which minimizes the function $g_s(x)$ on $\Zbold_+$.

Below and in the remainder of the paper, for notational simplicity,
we take $p=1/s$, for $s>0$. Let us now state the main theorems:

\begin{Theorem}[(Hopcount and weight asymptotics)]
\label{main_theorem}
For any fixed $s>0$ with $s\notin\SSS$, let $k^*(s)\in\{\lfloor
s+1\rfloor,\lceil s+1\rceil\}$ denote the unique integer that minimizes
the function defined in \eqref{def_g}. Then:
\begin{longlist}[(a)]
\item[(a)] the hopcount $H_n=H_n(s)$ converges in probability to $k^*(s)$ as
$n\to\infty$:
\[
\prob\bigl(H_n = k^*(s)\bigr) \to1;\vadjust{\goodbreak}
\]
\item[(b)] the optimal weight $W_n=W_n(s)$,
properly normalized converges in distribution as $n\to\infty$,
\begin{eqnarray*}
&&\prob\biggl(\frac{k-1}{sg_s(k)} (\log{n})^{s+1}\biggl(W_n - \frac{g_s(k)}
{(\log{n})^s}\biggr) + \frac{k-1}{2}\log{\log{n}} -\frac{p (k-1) \log
{g_s(k)}}{2} > t\biggr)\\
&& \quad\to\exp(-a_k \e^t), \qquad t\in\Rbold,
\end{eqnarray*}
where $k=k^*(s)$, and the sequence of constants
$(a_k)_{k\ge1}$ is defined by
%
%
\begin{equation}
\label{eqn:a-def}
a_{k}=\biggl(\frac{2\uppi p}{1+p}\biggr)^{(k-1)/2}k^{((k-1)p-1)/2}.
\end{equation}
\end{longlist}
\end{Theorem}

Theorem~\ref{main_theorem} states that the hopcount
$H_n$ converges to the optimal value of the function $x\mapsto g_s(x)$
defined in \eqref{def_g}, while the rescaled and recentered minimal
weight $W_n$ converges in distribution to a Gumbel distribution.
We can intuitively understand this as follows. For fixed~$k$,
the minimal path of length $k$ is similar to an
\textit{independent} minimum of copies of sums of $k$
random variables $E^{-s}$. The number of independent copies is equal
to the number of disjoint paths between vertices $1$ and $n$, which is
of order $n^{k-1}$. While
on $K_n$, the appearing paths do \textit{not} have independent
weights, the paths that are particularly short are almost independent.
Now, the independent problem can be handled in two steps. First, we analyze
the behavior of the random variable $Z_k=E_1^{-s}+\cdots+E_k^{-s}$.
In this analysis, the function $g_s$ appears in the lower tail
of the distribution.
Second, we study the asymptotics of the minimum of $n^{k-1}$
of such random variables, which can be seen to be of order
$g_s(k)/(\log{n})^s$. This explains why the minimal integer
value of $g_s$ is the crucial value for the hopcount, while
the minimum of a large number of independent random variables with
distribution $Z_k$,
properly rescaled and recentered, converges to a Gumbel distribution
by standard extreme value arguments. This intuitively explains
Theorem~\ref{main_theorem}. The main difficulty in the proof is to
handle the
fact that the weights of paths in the complete graph are actually
\textit{not}
independent, and we use the method of Stein--Chen for the Poisson
approximation to
deal with the available dependence.

Let us now deal with the case where $s\in\SSS$.
\begin{Theorem}[(The special set $\SSS$)]
\label{theo:main-s-in-ss}
Suppose $s \in\SSS$, so that both $\lfloor s+1 \rfloor$ and $\lceil
s+1 \rceil$ minimize
$g_s(\cdot)$ over $\Zbold_+$.
Define a sequence of independent random variables $(\Xi_k)_{k\geq2}$,
where, for any $k\geq2$, $\Xi_k$ has the Gumbel survival function
%
%
\begin{equation}
\label{eqn:gumb}
\prob(\Xi_k > t) = \exp(-a_k\e^t), \qquad t\in\Rbold,
\end{equation}
with $(a_k)_{k\geq2}$ defined in \eqref{eqn:a-def}. Then:
\begin{longlist}[(a)]
\item[(a)] the optimal weight $W_n=W_n(s)$, properly normalized, converges
in distribution as
\mbox{$n\rightarrow\infty$}. More precisely
\[
\frac{(\log{n})^{s+1}}{s g^*}\biggl(W_n - \frac{g^*}
{(\log{n})^s}\biggr) + \frac{1}{2}\log{\log{n}} -\frac{p \log{g^*}}{2}
\convd\min\biggl(\frac{\Xi_{\lfloor s+1 \rfloor}}{\lfloor s+1 \rfloor
-1},\frac{\Xi
_{\lceil s+1 \rceil}}{\lceil s+1 \rceil-1}
\biggr),
\]
where $g^*=g_s(\lfloor s+1 \rfloor)=g_s(\lceil s+1 \rceil)$;
\item[(b)] the hopcount $H_n=H_n(s)$ converges in distribution as
$n\rightarrow\infty$, that is
\[
H_n(s) \convd H^*,
\]
where
\[
H^* = \arg\min\bigl\{\Xi_k/(k-1)\dvt k\in\{\lfloor s+1 \rfloor, \lceil s+1
\rceil\}\bigr\}.
\]
\end{longlist}
\end{Theorem}

Another quantity of interest is the distribution of optimal paths between
vertex 1 and a \textit{set} of vertices.
In telecom, this is called multicast, since one source sends to a
multiple number of users.
This also follows from the analysis in the paper. We shall give a
brief idea of the proof in
Section~\ref{sec:multi-distance}. The result is stated for $s\notin
\SSS$,
but one could state an equivalent result for $s\in\SSS$ as well.
Before we state the result we need some further notation. Recall that
we used $k^*(s)$ to denote the unique minimizer of $g_s(\cdot)$ over
$\Zbold_+$. For fixed $m\geq1$, let $\{\eta_i\}_{1\leq i\leq m}$
denote independent copies of the Gumbel
random variable defined in \eqref{eqn:gumb} with $k=k^*(s)$.

\begin{Corollary}[(Multipoint distances)]
\label{corr:multi-point}
Fix $m\geq1$ distinct vertices say $2,3,\ldots, m+1$ in $K_n$. Suppose
$s\notin\SSS$, and let
$\{W^{(j)}_n\}_{2\leq j\leq m+1}$ denote the weight of the optimal
path from $1$ to these vertices. Write
\[
\tilde{W}^{(j)}_n = (\log{n})^{s+1}\biggl(W^{(j)}_n - \frac{g_s(k)}
{(\log{n})^s}\biggr) + \frac{k-1}{2}\log{\log{n}} -\frac{p (k-1) \log
{g_s(k)}}{2},
\]
where $k=k^*(s)$. Then, as $n\rightarrow\infty$,
\[
(\tilde{W}_n^{(j)})_{2\leq j\leq m+1}\convd(\eta_j)_{1\leq j\leq m}.
\]
\end{Corollary}

\textit{Organization of the paper}:
The paper is organized as follows. We first discuss the relevance of
our results and techniques in Section~\ref{sec:disc}. We shall then
continue to prove the main results in Section~\ref{sec:proofs}.

\section{Discussion}
\label{sec:disc}
We now provide a discussion of the various concepts used in this paper
and the relevance of the results.
\begin{longlist}[(a)]
\item[(a)] \textit{Stochastic mean-field model of distance}: This notion refers
to the complete graph with
exponentially distributed edge weights having unit mean. The model
gives a simpler but mathematically\vadjust{\goodbreak} more tractable model
of distances between random points in high dimensions.
While one can consider other edge distributions, the memoryless
property of exponential
random variable allows one to give clean proofs in a number of
different contexts,
including first passage percolation; see \cite{janson1999one} where
this property is used to great
effect to derive refined asymptotics.
We also refer to \cite{aldous2003objective}, where many other
computations are derived in this context
with the help of a powerful infinite random structure called the \textit{
weighted infinite tree}.
\item[(b)] \textit{Weak and strong disorder}: The last few years, with the
availability of an enormous amount of
data on real-world networks, has witnessed an explosion in network
models for these real-world
networks as well as dynamics on them. Physicists have been highly
interested in
understanding the effect of random disorder on the various flow
carrying properties
of these network models. Via simulations, they have predicted a number
of fascinating phenomena
in these networks. Regarding the notions of weak and strong disorder
mentioned in Section~\ref{sec-intro}, we refer the interested reader to
\cite{havlin2005optimal,braunstein2003optimal,braunstein2006optimal} and \cite{sreenivasan2004effect} and the
references therein.
\item[(c)] \textit{First passage percolation}: First passage percolation
problems have been of great interest to
probabilists for quite a while now, not just because of their origin
from physical motivations of modeling disordered random flow systems,
but also because this process and its variants (e.g., oriented first
passage percolation and last passage percolation) arise as basic
constructing blocks for more complicated problems, such as the contact
process. There has been an intensive study of this model on the
$d$-dimensional lattice (see, e.g., \cite{kesten1986aspects,smythe1978first} and \cite{howard2004models}). The case of the
complete graph with exponential edge weights was analyzed in \cite
{janson1999one}, where, in particular, it was proved that the weight
and hopcount of the optimal path satisfy
\[
nW_n-\log{n}\convd\Xi,
\]
and
\[
\frac{H_n}{\log{n}} \convp1,
\]
as $n\to\infty$, where $\Xi\in\Rbold$ is a non-degenerate random variable.

In the last few years, due to the connections to real-world networks
described above, these questions have taken on an added significance,
and a number of studies both at the non-rigorous level \cite
{braunstein2003optimal} and rigorous level (see, e.g., \cite
{bhamidi2010first}) have been undertaken to study such questions in
many other random graph models.
\item[(d)] \textit{Proof techniques}: A number of different techniques have
been used in the analysis of
first passage percolation asymptotics in various contexts, ranging from
subadditivity methods in the context
of the lattice, to continuous-time branching process embeddings and
renewal theory in the context of
various random graph models. The paper \cite{BH09} used embeddings
into a
particular continuous-time branching process to derive the results in
\eqref{eqn:hns-s-big-0} and \eqref{eqn:wn-s-big-0}.
As far as we know, the present paper is the first paper that uses the
method of Stein--Chen for Poisson approximation
to derive refined asymptotics in the first passage percolation context.
In particular the results here complete the program started in \cite{BH09}
and show that for $\beta\leq0$, $H_n(\beta)$ is a tight sequence of
random variables,
while for $\beta>0$, $H_n(\beta)/(\beta\log{n})\convp1$.
Clearly, this further shows that there are at least \textit{two} universality
classes for first-passage percolation on the complete graph in terms of
the edge weight distribution.
When $\beta\leq0$, the weights $E^{\beta}$ are in the same
universality class
as the weight 1, in the sense that the hopcount remains bounded,
while for $\beta>0$, they are in the same universality class as the
exponential distribution
arising for $\beta=1$. As discussed in more detail in \cite{BH09}, this
raises the question of
what the universality classes for first passage percolation on $K_n$ are.
In particular, does $H_n$ always satisfy a central limit theorem
whenever $H_n\rightarrow\infty$? Or are there classes of edge weight
distributions
where the behavior is even different? For example, is there a class
of random edge weights where the behavior is similar as for the
minimal spanning tree, where $H_n$ is of the order $n^{1/3}$.
\item[(e)] \textit{Multi-point distances and exchangeability}: The classical
probability theory of
exchangeability has been used in the last few years to analyze various
complex random structures; see \cite{aldous2009more} for a nice,
modern survey.
In the context of Corollary~\ref{corr:multi-point}, one can analyze
such questions in a number of different contexts (such as the
stochastic mean-field model of distance or your favorite random graph
model with your favorite random edge weights). For the stochastic
mean-field model,
the multi-point optimal path weights converge (after proper rescaling
and recentering) to an \textit{exchangeable} sequence of random variables.
In the present model, we can once again show convergence but to an
independent sequence of random variables.
\end{longlist}

\section{Proofs}
\label{sec:proofs}
This section contains the proofs of the main results.
We start with an outline of the proof. In Section~\ref{sec:tight-hop},
we shall show that the hopcount $H_n(s), s>0,$
is a tight sequence
of random variables as $n\to\infty$. In Section~\ref{minima},
we shall derive the asymptotic behavior, for $z \downarrow0$, of
the distribution function $F_k(z)$, where
%
\begin{equation} \label{dfk} F_k(z)=\prob(E_1^{-s}+E_2^{-s}+\cdots
+E_k^{-s}\le z),
\end{equation}
and where $E_1,E_2,\ldots$ is an i.i.d.~sequence of $\EXP(1)$ random
variables.
Denoting the weight of the minimal path with exactly $k$ edges between
vertex $1$ and vertex $n$
by $W_k(n)$, we then show that for each $K$, and each $0<\vep<1$,
uniformly in $2\leq k\le K$, w.h.p.,
%
\begin{equation} \label{firstlower} (\log n)^s W_k(n)\geq(1-\vep)
g_s(k),
\end{equation}
where $g_s(x)$ is defined in \eqref{def_g}.

Inequality \eqref{firstlower} yields a first-order lower bound for all
values of $k\ge2$. We will show in the paper that the function
$g_s(x)$ determines the first-order asymptotics of the weights
$W_k(n)$. The behavior of $g_s$ near its minimum value
determines the asymptotic behavior of the hopcount $H_n$. Roughly
speaking, the hopcount $H_n$
will converge in probability to the integer $k=k^*(s)$
that minimizes the function $g_s(x), x\ge2$, over the set $\Zbold
_+$. The above statement about the convergence of $H_n$ is true for
every $s>0$ for which the minimizing integer of $g_s(x)$ for $x\geq2$
is unique,
that is, $s\notin\SSS$. In Section~\ref{sec-asym}, we study the
minimum of
an independent number of $n^{k-1}$ random variables. Each of these
$n^{k-1}$ variables
is the sum of $k$ i.i.d.~random variables with distribution
$E^{-s}$. The result is used to complete the proof of Theorem \ref
{main_theorem}
when $0<s\leq1$. In Section~\ref{sec_greater1}, we extend the analysis
to $s>1$ and complete the proof of Theorem~\ref{main_theorem}
by studying the second order asymptotics of the minimal weight of
paths of length $k$ in the complete graph $K_n$.

For $s_j\in\SSS$, to decide whether the hopcount $H_n$ converges in
probability either to $\lfloor s_j+1\rfloor$
or to $\lceil s_j+1\rceil$, we need the second order asymptotics of $W_k(n)$,
which is carried out in detail in Section~\ref{sec_bad}.
In Section~\ref{sec:multi-distance}, we sketch the proof of Corollary
\ref{corr:multi-point}.

\subsection{Tightness of the hopcount}
\label{sec:tight-hop}
Note that the minimal weight $W_n$ satisfies the following inequality:
%
\begin{equation} \label{lowerbndWn} W_n\ge H_n \cdot\min_{1\le j \le
n(n-1)/2} E_j^{-s},
\end{equation}
where $E_j\sim\EXP(1)$ are independent.
Since the maximum of $n$ independent exponentials scales like $(1+\mathrm{o}_{
\prob}
(1))\log n$, we obtain
from (\ref{lowerbndWn}) that w.h.p.
%
\begin{equation} \label{lowerbndWn2} W_n\ge H_n \log^{-s}
\bigl(n(n-1)/2\bigr)\bigl(1+\mathrm{o}_{ \prob}(1)\bigr).
\end{equation}
On the other hand, $W_n$ is, at most, equal to the minimal weight of
all two-edge paths between $1$ and $n$. Here, a two-edge path is a path
of the form
$1\to j \to n, j=2,3,\ldots, n-1$, so that
%
\begin{equation} \label{upperWn} W_n \leq\min_{2\le j\le n-1}
\bigl((E'_j)^{-s}+(E''_j)^{-s}\bigr),
\end{equation}
where $E'_j, 2\le j \le n-1,$ and $E''_j, 2\le j \le n-1,$ are
independent $\EXP(1)$ random variables.
It is not hard to verify (see Lemma~\ref{dist+sums} in Section \ref
{minima}) that
\eqref{upperWn} implies that w.h.p.,
%
\begin{equation} \label{upperbndWn2} W_n\le\frac{C}{(\log n)^{s}}.
\end{equation}
Inequalities \eqref{lowerbndWn2} and \eqref{upperbndWn2} together imply
that w.h.p.,
%
\begin{equation} \label{tightnessHn} H_n\leq C (\log2)^s.
\end{equation}
We conclude that $H_n$ is a tight sequence of random variables.
The remainder of the proof will reveal that, in fact, $H_n$ converges in
distribution, either to a constant $k^*(s)$ when $s\notin\SSS$,
or to a random variable giving positive mass to two values when
$s\in\SSS$.

\subsection{The first-order lower bound}
\label{minima}
We start with an investigation of the distribution function $F_k$ of an
\textit{independent} sum of $k$ inverse powers of exponentials, that is,
%
\begin{equation} \label{def_Zk} Z_k=E_1^{-s}+\cdots+E_k^{-s},\qquad s>0.
\end{equation}
\begin{Lemma}[(Sums of inverse powers of exponentials)]
\label{dist+sums}
Fix $s>0$, and put $p=1/s$. Then, for $z\downarrow0$,
%
\begin{equation} \label{behav_Fk} F_k(z)\sim a_k z^{-(k-1)p/2} \e^{-
k^{p+1} z^{-p}},
\end{equation}
where $(a_k)_{k\geq1}$ is defined in \eqref{eqn:a-def},
and where, for arbitrary real functions $g$ and $h$, $g(z)\sim h(z),
z \downarrow0$,
means that $\lim_{z\downarrow0}g(z)/h(z)=1$.
\end{Lemma}

\begin{pf}The result for $k=1$ is immediate from $F_1(z)=\e
^{-z^{-p}}, z>0$.
We proceed by induction. Suppose that \eqref{behav_Fk} holds for some
integer $k\ge1$, then
\begin{eqnarray*}
F_{k+1}(z)&\sim& \int_0^z a_k
(z-y)^{-(k-1)p/2}
\e^{- k^{p+1} (z-y)^{-p}}
\,\mathrm{d}(\e^{-y^{-p}})\\
&=&pa_kz^{-(k-1)p/2}z^{-p}\int_0^1 x^{-p-1}(1-x)^{-(k-1)p/2}\e
^{-z^{-p}h_k(x)} \,\mathrm{d}x,
\end{eqnarray*}
where $h_k(x)=x^{-p}+k^{p+1}(1-x)^{-p}$. The function $h_k$ has a
minimum at $x=1/(k+1)$, since
$h_k(1/(k+1))=(k+1)^{p+1}$, $h'_k(1/(k+1))=0,$ and
\[
h''_k\bigl(1/(k+1)\bigr)=p(p+1)(k+1)^{p+2}\biggl(1+\frac{1}{k}\biggr).
\]
Hence, from a standard Laplace-method argument, we obtain
\begin{eqnarray*}
 F_{k+1}(z)&\sim& pa_kz^{-(k-1)p/2}z^{-p}\int_0^1
x^{-p-1}(1-x)^{-(k-1)p/2}\e^{-z^{-p}h_k(x)} \,\mathrm{d}x\\
& \sim&pa_kz^{-(k-1)p/2}z^{-p}
(k+1)^{p+1}\bigl(k/(k+1)\bigr)^{-(k-1)p/2}\\
&&{}\times  \e^{-z^{-p}h_k(1/(k+1))}
\sqrt{\frac{2\uppi}{z^{-p}h_k''(1/(k+1))}}.
\end{eqnarray*}
From the latter expression, we obtain
\[
a_{k+1}=\biggl(\frac{2\uppi p}{1+p}\biggr)^{1/2}
a_k
(k+1)^{(kp-1)/2}
k^{-((k-1)p-1)/2}.
\]
This recursion is telescoping in $k$. Defining $c=(\frac{2\uppi
p}{1+p})^{1/2}$
and $b_k=k^{((k-1)p-1)/2}$, we find
\[
a_{k+1}=c a_k \frac{b_{k+1}}{b_k}\quad\Rightarrow\quad a_{k+1}=c^k \frac
{b_{k+1}}{b_1}a_1=c^k b_{k+1},
\]
which yields (\ref{eqn:a-def}).
\end{pf}

Using the above lemma, we obtain the following first-order lower bound
for $W_k(n)$:
\begin{Theorem}[(First-order lower bound)]
\label{lowerbound}
Fix $s>0$ and an arbitrary large integer $K$. For each $0<\vep<1$,
with the
function $g_s$ defined in \eqref{def_g}, w.h.p. and uniformly
in $k
\in\{2,3, \ldots,K\}$,
\[
(\log n)^s W_k(n)\ge(1-\vep)g_s(k).
\]
\end{Theorem}

\begin{pf} Fix $s>0$ and $2\le k < n$, and define, for $0<\vep<1$,
%
\begin{equation} x_{k,n}=x_{k,n}(\vep)=(1-\vep) \frac{g_s(k)}{(\log
n)^s}.
\end{equation}
Let $N_k^{(n)}(x), x>0,$ be the number of paths between $1$ and
$n$ with exactly $k$ edges and weight at most $x$.
Note that the total number of paths with exactly $k$ edges between $1$
and $n$ is $\prod_{j=2}^k (n-j)$. Thus, according to Lemma \ref
{dist+sums}, for $x\downarrow0$,
%
\begin{equation} \label{exp_good}
\expec\bigl[N_k^{(n)}(x)\bigr]=\Biggl[\prod
_{j=2}^k (n-j)\Biggr] F_k(x)\sim \Biggl[\prod_{j=2}^k (n-j)\Biggr] a_k x^{-(k-1)p/2} \e
^{- k^{p+1} x^{-p}}.
\end{equation}
For $n\to\infty$ the expression $x_{k,n}\downarrow0,$ and the term
$x_{k,n}^{-(k-1)p/2}$ blows up only polynomially fast,
while $\exp\{- k^{p+1} x_{k,n}^{-p}\}$ tends to $0$ exponentially fast.
Using that $\prod_{j=2}^k (n-j)< n^{k-1}$ and abbreviating $N_k^{
(n)}=N_k^{(n)}(x_{k,n})$,
we reach to the conclusion that
\[
\expec\bigl[N_k^{(n)}\bigr]\leq n^{k-1}\exp\{- k^{p+1} x_{k,n}^{-p}\}=\exp
\biggl\{
-\biggl(\frac1{(1-\vep)^p}-1\biggr)(k-1)\log n
\biggr\}.
\]
Boole's inequality and the Markov inequality together yield
\begin{eqnarray*}
&&\prob\Biggl(\bigcup_{k=2}^{K}\{(\log n)^s
W_k(n)<(1-\vep)g_s(k)\}\Biggr)\\
&&\quad\le\sum_{k=2}^{K}\prob\bigl((\log n)^s W_k(n)<(1-\vep)g_s(k)\bigr)
 \le\sum_{k=2}^{K}\prob\bigl(N_k^{(n)}\ge1\bigr)\\
&&\quad\le\sum_{k=2}^{K}\expec\bigl[N_k^{(n)}\bigr] \le\sum_{k=1}^{\infty}\exp
\biggl\{-\biggl(\frac1{(1-\vep)^p}-1\biggr)k\log n\biggr\}.
\end{eqnarray*}
Since the summand on the right-hand side is of order $n^{-p\vep k}$, we
may conclude that the probability that $(\log n)^s W_k(n)<(1-\vep
)g_s(k)$, for some $2\le k\le K$,
tends to $0$ as $n\to\infty$.
\end{pf}

\subsection{Second-order asymptotics}
\label{sec-asym}
In this section we identify the second-order asymptotics
for the \textit{independent} minimum of $n^{k-1}$ random variables, where
each of
these random variables has distribution function $F_k(z), z>0$.
The proof of Theorem~\ref{main_theorem} for $0<s \le1$
follows quite easily from this and the lower
bound \eqref{firstlower}. The proof of Theorem~\ref{main_theorem}
for $s>1$ is postponed to the next section.

We write
%
\begin{equation} \label{indep_wk} W_k^{(\mathrm{ind})}=\min_{1\le j
\le n^{k-1}} Y_{k,j},
\end{equation}
where $Y_{k,1},\ldots,Y_{k,n^{k-1}}$ are i.i.d. with distribution
function $F_k$. The following theorem derives the
asymptotics of $W_k^{(\mathrm{ind})}$:

\begin{Theorem}[(Minimum for independent r.v.s)]
\label{min_independent}
For $k\ge2$ fixed, the minimal weight $W_k^{(\mathrm{ind})}$ defined
in \eqref{indep_wk}, satisfies
%
\begin{eqnarray} \label{weightas}
&& \prob\biggl(\frac{k-1}{sg_s(k)}(\log
n)^{s+1} \biggl( W_k^{(\mathrm{ind})}-\frac{g_s(k)}{(\log n)^s}\biggr) +\frac
{k-1}{2}\log{\log n} -\frac{p(k-1) \log{g_s(k)}}{2}>t\biggr)\qquad\quad
\nonumber
\\[-8pt]
\\[-8pt]
\nonumber
&&\quad \to\e
^{-a_k\e^t},
\end{eqnarray}
where $(a_k)_{k\geq1}$ is defined in \eqref{eqn:a-def}.
\end{Theorem}

\begin{pf} We compute $z_n=z_n(t)$, such that
\[
\bigl(1-F_k(z_n)\bigr)^{n^{k-1}}\to\exp\{-a_k\e^t\}.
\]
Taking logarithms on both sides and using that $\log\{1-F_k(z_n)\}\sim
-F_k(z_n)$
for $z_n\rightarrow0$, this is equivalent to
%
\begin{equation} \label{asymFkzn} n^{k-1}F_k(z_n)\to a_k\e^t,
\end{equation}
or
%
\begin{equation} \label{gerard1} (k-1)\log{n}+\log\{F_k(z_n)\}\to
t+\log a_k.
\end{equation}
Put $z_n=\kappa(\log{n})^{-s}+\zeta_n(t)$, where $\kappa=g_s(k)$.
From Lemma~\ref{dist+sums}, we find that \eqref{gerard1} is
equivalent to
%
\begin{eqnarray} \label{gerard2}
&&(k-1)\log{n} -(k-1)p/2 \log\bigl(\kappa
(\log{n})^{-s}+\zeta_n(t)\bigr)
\nonumber
\\[-8pt]
\\[-8pt]
\nonumber
&&\quad{}-k^{p+1}\bigl(\kappa(\log{n})^{-s}+\zeta
_n(t)\bigr)^{-p}\to t.
\end{eqnarray}
Writing
\[
\kappa(\log{n})^{-s}+\zeta_n(t)
=
\kappa(\log{n})^{-s}\bigl(1+\zeta_n(t)(\log{n})^{s}/\kappa\bigr),
\]
yields
\begin{eqnarray*}
&&(k-1)\log{n}
-(k-1)p/2 \log\bigl(\kappa(\log{n})^{-s}\bigl(1+\zeta_n(t)(\log
{n})^{s}/\kappa\bigr)\bigr)\\
&&\quad{} -k^{p+1}\kappa^{-p}\log{n}\cdot\bigl(1+\zeta_n(t)(\log
{n})^{s}/\kappa\bigr)^{-p}\to t.
\end{eqnarray*}
Using that $k^{p+1}\kappa^{-p}=k-1$ and $ps=1$, we arrive at
\begin{eqnarray*}
&&
(k-1)\log{n}+(k-1)/2 \log(\log{n})
-(k-1)p/2\log
\bigl(1+\zeta_n(t)(\log{n})^{s}/\kappa\bigr)
\\
&&\quad{}
-(k-1)\log{n}\cdot\bigl(1+\zeta_n(t)(\log{n})^{s}/\kappa
\bigr)^{-p}\to t
+(k-1)p/2 \log\kappa.
\end{eqnarray*}
Now we choose
%
\begin{equation} \zeta_n(t)= (\log n)^{-s-1}\cdot\bigl(\zeta\log{\log
{n}}+h(t)\bigr)\quad \mbox{or} \quad\zeta_n(t)(\log n)^s=\frac{(\zeta\log{\log
{n}}+h(t))}{\log n}.\ \
\end{equation}
Then
\[
(k-1)p/2\log
\bigl(1+\zeta_n(t)(\log{n})^{s}/\kappa\bigr)
=\mathrm{O}\biggl(\frac{\log{\log n}}{\log n}\biggr)\to0,
\]
and
\begin{eqnarray*}
-(k-1)\log{n}\cdot\bigl(1+\zeta_n(t)(\log{n})^{s}/\kappa
\bigr)^{-p}&\sim&
-(k-1)\log{n}\cdot\bigl(1-p\zeta_n(t)(\log{n})^{s}/\kappa\bigr)\\
& = &-(k-1)\log{n}+\frac{(k-1)p}{\kappa}\bigl(\zeta\log{\log{n}+h(t)}\bigr),
\end{eqnarray*}
resulting in
%
\begin{equation} \label{gerard3} \zeta=-\kappa/2p \quad\mbox{and}\quad \frac
{(k-1)ph(t)}{\kappa}=t+\frac12 (k-1)p \log\kappa.
\end{equation}
Hence, we can choose
%
%
\begin{eqnarray}
\label{corr_zn}
z_n(t)&=&g_s(k)(\log{n})^{-s}+\zeta_n(t)\nonumber\\
&=&
g_s(k)(\log{n})^{-s}+(\log n)^{-s-1}\cdot\bigl(\zeta\log{\log
{n}}+h(t)\bigr)\\
&=&
\frac{g_s(k)}{(\log{n})^s}
+ \frac{g_s(k)}{(\log{n})^{s+1}}\biggl[ - \frac{\log{\log{n}}}{2p}+
\frac{t}{(k-1)p} +
\frac{\log{g_s(k)}}{2}\biggr].\nonumber
\end{eqnarray}
\upqed\end{pf}

We now turn to the proof of Theorem~\ref{main_theorem} in the
case where $0<s\le1$:

\begin{pf*}{Proof of Theorem \protect\ref{main_theorem} in case $0<s\le1$}
Observe from
Theorem~\ref{lowerbound} that for any $K$, w.h.p. and
uniform in $k\in\{2,3,\ldots,K\}$,
%
\begin{equation} \label{ger1} (\log n)^s W_k(n)\ge(1-\vep
)g_s(3)>g_s(2),
\end{equation}
where the latter inequality follows since for the indicated values of
$s$, the function
$g_s$ is increasing on $[2,\infty)$ and where we can take
$\vep<\min_{0<s\le1} [1-g_s(2)/g_s(3)]= 1/9$.
On the complete graph with $n$ vertices the paths of length
$2$ have \textit{independent} total weight, since they are disjoint.
The number of paths of length 2 is equal to $n-2\sim n$,
so that we can conclude from the Theorem~\ref{min_independent} that,
for any $\vep>0$ and w.h.p.,
%
\begin{equation} \label{ger2} (\log n)^sW_2(n)\in\bigl(g_s(2)-\vep
,g_s(2)+\vep\bigr).
\end{equation}
From \eqref{ger1} and \eqref{ger2} it is immediate that,
{w.h.p.}, the
minimal-weight path
is either a path of length $2$ or has a length exceeding $K$.
Since $K$ can be taken arbitrary large and the hopcount $H_n$ is tight
(see~\ref{tightnessHn}),
we conclude that
\[
H_n(s)\convp2
\]
for $0<s\le1$. Consequently, $W_n=W_n(2)$, w.h.p.,
and statement (b) of Theorem~\ref{main_theorem} follows from \eqref
{weightas} for $0<s\le1$ and $k=2$.
\end{pf*}

\subsection{The case $s>1$}
\label{sec_greater1}
In this section we treat the case $s>1$.
The number of paths with $k\geq2$ edges between the
vertices~$1$ and $n$ is equal to $\prod_{j=2}^k (n-j) \sim n^{k-1}$.
Let $\SSS_k(n)$ denote the set of all such paths.\vspace*{1pt}
As before, we let $F_k$ denote the distribution function
of the sum of $k$ independent random variables each with distribution
equal to
the distribution of $E^{-s}$, and by $ N^{(n)}_k(z), z>0,$ the
number of paths with $k$ edges which
have total weight $w_s(\PP)=\sum_{e\in\PP} E_e^{-s}$ less than~$z$.
Recall the definition of $z_n(t)$ in \eqref{corr_zn}.
From Theorem~\ref{min_independent} and its proof (compare \eqref{asymFkzn}),
we conclude that, as $n\rightarrow\infty$,
%
%
\begin{equation}
\label{eqn:lambdak-asymp}
\lambda_k^{(n)}(t)
:=
\expec\bigl[
N_k^{(n)}
(z_n(t))\bigr]
\sim
n^{k-1} F_k(z_n(t))
\to a_k \e^t.
\end{equation}
We shall prove the following proposition:
\begin{Proposition}[(Poisson approximation for small weight paths)]
\label{prop:poisson-dist-err}
Fix $s>1$ and let $\Poi_k^{(n)}(t)$ be a Poisson random variable
with mean $\lambda_k^{(n)}(t)$. Then,
both for $k=\lfloor s+1 \rfloor$ and $k=\lceil s+1 \rceil$,
as $n\rightarrow\infty$,
\[
d_{\mathrm{TV}}\bigl(N_k^{(n)}(z_n(t)) , \Poi_k^{(n)}(t)\bigr) \to0,
\]
where $d_{\mathrm{TV}}$ denotes the total variation distance.
\end{Proposition}

Assuming the proposition let us first show how to complete the proof of
Theorem~\ref{main_theorem}.

\begin{pf*}{Proof of Theorem \protect\ref{main_theorem} in case $s>1$ and
$s\notin
\mathcal{S}$}
Observe that
%
\begin{equation} \label{duality}
\prob\bigl(W_k(n)> z_n(t)\bigr) = \prob
\bigl(N_k^{(n)}(z_n(t)) = 0\bigr).
\end{equation}
Now Proposition~\ref{prop:poisson-dist-err}
together with \eqref{eqn:lambdak-asymp} imply that, for $k=\lfloor s+1
\rfloor$ and $k=\lceil s+1 \rceil$, as $n\to\infty$,
%
\begin{equation} \label{second_order_as}
\prob\bigl(W_k(n)> z_n(t)\bigr)\to
\exp(-a_k \e^t).
\end{equation}
Note that the weak convergence shows in particular that $(\log
n)^sW_k(n)$ converges in probability to $g_s(k)$ for the two indicated
values of $k$.
This together with the lower bound proven in Theorem~\ref{lowerbound},
and an argument similar to the case $0<s\le1$ then completes the proof
of Theorem~\ref{main_theorem}, in case the integer that minimizes
$g_s(x)$ is unique, that is, in case
$s\notin\mathcal{S}$.
\end{pf*}

\begin{pf*}{Proof of Proposition \protect\ref{prop:poisson-dist-err}}
We shall use \cite{BHJ92}, Theorem 1.A. Before quoting this result, we
shall need to setup some notation.
Let $\mathcal{I}$ be a finite index set, and let $\{I_\alpha\dvtx \alpha\in
\mathcal{I}\}$
be a family of indicator random variables and write $p_\alpha= \expec
[I_\alpha]$. Let
\[
W=\sum_{\alpha\in\mathcal{I}} I_\alpha,\qquad \lambda= \expec[W] = \sum
_{\alpha\in\mathcal{I}} p_\alpha.
\]
Now suppose for each $\alpha$ we can decompose the index set
$\mathcal{I}$ as $\mathcal{I}= \{\alpha\}\,\cup\,\IIst(\alpha) \,\cup\,\SSst
(\alpha)$,
where we shall think of $\{I_\beta\dvt \beta\in\IIst(\alpha)\}$ to be
the set of random variables
which ``strongly depend'' on~$I_\alpha$ while $\{I_{\beta^\prime
}\dvt \beta^\prime\in\SSst(\alpha)\}$
consists of the set of random variables which only ``weakly depend'' on~$I_\alpha$.
Now let $Z_\alpha= \sum_{\beta\in\IIst(\alpha)} I_\beta$, while
\[
Y_\alpha= W-I_\alpha- Z_\alpha= \sum_{\beta^\prime\in\SSst
(\alpha
)} I_{\beta^\prime}.
\]
Then with this notation,
the following is just one example of the power of the Stein--Chen
machinery for
Poisson approximation for weakly dependent indicator random variables:
\begin{Theorem}[(Stein--Chen Poisson approximation (\cite{BHJ92}, Theorem 1.A))]
\label{theo:stein-chen-poi}
With the above notation,
\[
d_{\mathrm{TV}}(W,\Poi(\lambda))
\leq
\min(1,\lambda^{-1})
\sum_{\alpha\in\mathcal{I}}
(
p_\alpha^2 + p_\alpha\expec[Z_\alpha]+ \expec[I_\alpha Z_\alpha]
)
+\min(1,\lambda^{-1/2})
\sum_{\alpha\in\mathcal{I}} \eta_\alpha,
\]
where $\eta_\alpha$ is such that
\[
|\expec[I_\alpha g(Y_\alpha+1)] - p_\alpha\expec[g(Y_\alpha
+1)]|\leq
\eta_\alpha\Vert g\Vert ,
\qquad
\alpha\in\mathcal{I}
\]
for all bounded functions $g$ on $\Zbold_+$, and where $\Vert \cdot\Vert $ is
the supremum norm.
\end{Theorem}

To apply Theorem~\ref{theo:stein-chen-poi} to the situation at hand,
we take
$\mathcal{I}=\SSS_k(n)$, the set of paths between $1$ and $n$ having
precisely $k$ edges.
For $\alpha\in\SSS_k(n)$, we denote by
%
\begin{equation} \label{I-alpha-def} I_{\alpha}=I_\alpha(z_n(t))=
\mathbh{1}_{\{w_s(\alpha)\leq z_n(t)\}},
\end{equation}
where, as before, $w_s(\alpha)=\sum_{e\in\alpha} E_e^{-s}$
denotes the weight of the path $\alpha$,
and where $\mathbh{1}_{A}$ denotes the indicator of event $A$.
Furthermore,
let $p_k^{(n)}(t)$ denote
the expectation of $I_\alpha(z_n(t))$, that is,
%
\begin{equation} \label{pkn-def} p_k^{(n)}(t)=\prob\bigl(w_s(\alpha)\leq
z_n(t)\bigr)=F_k(z_n(t)).
\end{equation}
Let $\IIst(\alpha)\subseteq\SSS_k(n)$ denote the set of paths (not
including $\alpha$)
which have at least one edge in common with $\alpha$ (\textit{i.e.,
$\IIst
(\alpha)$ is the set
of paths $\beta$ for which
$I_\beta$ is ``strongly'' dependent on $I_\alpha$}), and let
$\SSst(\alpha)\subseteq\SSS_k(n)$ denote the set of paths that do not
overlap on any edge with $\alpha$.\vadjust{\goodbreak}
Note that the random variable $w_s(\alpha)$ is independent of
$\{w_s(\beta)\dvt\beta\in\SSst(\alpha)\}$.
Finally, in the above notation, note that
\[
Z_\alpha= \sum_{\beta\in\IIst(\alpha)}\mathbh{1}_{\{
w_s(\beta)\leq z_n(t)\}}.
\]
The independence of $w_s(\alpha)$ and $\{w_s(\beta)\dvt\beta\in
\SSst(\alpha)\}$
implies that we can take $\eta_\alpha= 0$ in Theorem~\ref
{theo:stein-chen-poi}.
Thus applying Theorem~\ref{theo:stein-chen-poi}, we get
%
%
\begin{eqnarray}
\label{eqn:dtv-bound1}
d_{\mathrm{TV}}\bigl(N_k^{(n)}(z_n(t)) , \Poi_k^{(n)}(t)\bigr) &\leq&
\frac{\sum_{\alpha\in\SSS_k(n)}[(p_k^{(n)}(t))^2+p_k^{
(n)}(t)\expec[Z_\alpha]+\expec[I_\alpha Z_\alpha]
]}{\lambda_k^{(n)}(t)}
\nonumber
\\[-8pt]
\\[-8pt]
\nonumber
&=& p_k^{(n)}(t)+\expec[Z_\alpha] + \frac{\expec[I_\alpha Z_\alpha
]}{p_k^{(n)}(t)},
\end{eqnarray}
where the last equality follows since $\lambda_k^{(n)}(t)=|\SSS
_k(n)|p_k^{(n)}(t)$
and since $\expec[Z_\alpha]$ and $\expec[I_\alpha Z_\alpha]$ are
independent of $\alpha$.
As before, by the choice of $z_n(t)$,
\[
n^{k-1} p_{k}^{(n)}(t)\to a_k \e^{t}.
\]
Thus, in particular, $p_{k}^{(n)}(t)\to0,$ as $n\to\infty$. Further,
there exists a constant $C_k$ such that, as $n\rightarrow\infty$,
\[
\expec[Z_\alpha]= |\IIst(\alpha)|p_{k}^{(n)}(t)\leq C_k n^{k-2}
p_{k}^{(n)}(t) \to0.
\]
Thus, the first two terms in \eqref{eqn:dtv-bound1}
vanish as $n\to\infty$. The last term requires some more analysis.
We note that
%
%
\begin{equation}
\expec[I_\alpha Z_\alpha]=\sum_{j=1}^{k-2} |\IIst_{k,j}(\alpha
)|p_{k,j}^{(n)}(t).
\label{eqn:st-sum-alpha}
\end{equation}
Here $\IIst_{k,j}(\alpha)\subseteq\SSS_k(n)$ consists of the set of
paths of length $k$ which overlap with $\alpha$
in exactly $j$ edges, while
\[
p_{k,j}^{(n)}(t)=
\prob\bigl(X_{k,k}\leq z_n(t), X_{k,j} \leq z_n(t)\bigr),
\]
where $X_{k,k}=\sum_{r=1}^k E_r^{-s},$ while
$X_{k,j}=\sum_{r=1}^j E_r^{-s}+ \sum_{r=j+1}^k \tilde{E}_r^{-s},
1\le
j \le k-2$
and $(E_i)_{i=1}^k$ and $(\tilde{E}_r)_{r=1}^k$ are two independent
vectors of
i.i.d. $\EXP(1)$ random variables.
We bound the probability $p_{k,j}^{(n)}(t)$ in the same way as
before, using the standard Laplace's method:
\begin{Lemma}[(Correlated sums of inverse powers of exponentials)]
\label{lemma:bi-bound}
Fix $k\geq3$, and let $1\leq i\leq k-2$. Then, for $z\downarrow0$,
\[
p_{k,j}^{(n)}(t)=\prob(X_{k,k}\le z, X_{k,i}\le z) \sim C_{k,i}
\frac{1}{z^{(k-i-1)p+ip/2}}
\exp\bigl(-z^{-p}[(k-i)\nu+i]^{p+1} \bigr),
\]
where $\nu=2^{1/(p+1)}$ and $C_{k,i}> 0$ is a constant.\vadjust{\goodbreak}
\end{Lemma}
\begin{pf}
The proof is given by straightforward computation using Laplace's method:
\begin{eqnarray*}
&&\prob(X_{k,k}\leq z,X_{k,i}\leq z)\\
&&\quad=
\prob\Biggl(\sum_{r=1}^k E_r^{-s}\leq z,\sum_{r=1}^i E_r^{-s}+\sum
_{r=i+1}^k ({\tilde E}_r)^{-s}\leq z\Biggr)=
\int_0^z
F^2_{k-i}(z-y) \,\mathrm{d}F_i(y)\\
&&\quad \sim a_ia_{k-i}^2\int_0^z
(z-y)^{-(k-i-1)p}
\e^{- 2(k-i)^{p+1} (z-y)^{-p}} \,\mathrm{d}y^{-(i-1)p/2}
\e^{- i^{p+1} y^{-p}}\\
&&\quad =a_ia_{k-i}^2\int_0^z
y^{-(i-1)p/2-1}\e^{- i^{p+1}
y^{-p}}\bigl(pi^{p+1}y^{-p}-(i-1)p/2\bigr)\\
&&\hspace*{42pt}\qquad{}\times (z-y)^{-(k-i-1)p}
\e^{- 2(k-i)^{p+1} (z-y)^{-p}} \,\mathrm{d}y
\\
&&\quad =a_ia_{k-i}^2 z^{-(k-i-1)p-(i-1)p/2}\int_0^1
x^{-(i-1)p/2-1}\bigl(pi^{p+1}x^{-p}z^{-p}-(i-1)p/2\bigr)\\
&&\hspace*{121pt}\qquad{} \times
(1-x)^{-(k-i-1)p}
\exp\{-z^{-p}h_{k,i}(x)\} \,\mathrm{d}x,
\end{eqnarray*}
where we abbreviate
\[
h_{k,i}(x)=i^{p+1}x^{-p}+2(k-i)^{p+1}(1-x)^{-p}.
\]
Put $\nu=2^{1/(p+1)}$. Then the minimum arises
in the point $x_{k,i}$, satisfying
$h_{k,i}^{\prime}(x_{k,i})=0$, which yields
\[
x_{k,i}=\frac{i}{(k-i)\nu+i}.
\]
Furthermore, $h_{k,i}(x_{k,i})=((k-i)\nu+i)^{p+1}$, while
\[
h''_{k,i}(x_{k,i})=\frac{p(p+1)}{i(k-i)\nu}\bigl((k-i)\nu+i\bigr)^{p+3}.
\]
Applying Laplace's method then yields
%
%
\begin{eqnarray}\label{zadelpunt}
&&\prob(X_{k,k}\leq z,X_{k,i}\leq z)\nonumber\\
&&\quad\sim
a_ia_{k-i}^2z^{-(k-i-1)p-(i-1)p/2}
(x_{k,i})^{-(i-1)p/2-1}\bigl(pi^{p+1}(zx_{k,i})^{-p}-(i-1)p/2
\bigr)\qquad\\
&&\qquad{} \times(1-x_{k,i})^{-(k-i-1)p}
\exp\{-z^{-p}h_{k,i}(x_{k,i})\}\sqrt{\frac{2\uppi}{z^{-p}h''_{k,i}(x_{k,i})}}.\nonumber
\end{eqnarray}
\upqed\end{pf}

Recall \eqref{eqn:dtv-bound1}. The first two terms on the right-hand
side vanish, as $n\to\infty$, hence
it suffices to show that
\[
\frac{\expec[I_\alpha Z_\alpha]}{p_k^{(n)}(t)} \to0.
\]
Using that $p_{k}^{(n)}(t)= \mathrm{O}(n^{-(k-1)})$
and that $|\IIst_{k,j}(\alpha)|\sim n^{k-j-1}$ it follows from \eqref
{eqn:st-sum-alpha} that
we now need to show for $1\leq j\leq k-2$,
\[
n^{2k-j-2} p_{k,j}^{(n)}(t) \to0,
\]
as $n\to\infty$. Now the polynomial terms ($z^k$ type terms) in
the approximation of $p_{k,j}^{(n)}(t) $ should not play a role.
Thus, using the fact that up to the first-order
%
\begin{equation} \label{zn-asympt} (z_n(t))^{-p}\sim\biggl(\frac
{g_s(k)}{\log^s{n}}\biggr)^{-p}=\frac{(k-1) \log n}{k^{1+p}},
\end{equation}
we need to show that for $1\leq j\leq k-2$ and with $\nu=2^{1/(p+1)},$
%
\begin{equation} \label{condPoiss} \biggl[\biggl(\biggl(1-\frac{j}{k}\biggr)\nu+ \frac
{j}{k}\biggr)^{p+1}- \biggl(2-\frac{j}{k-1}\biggr)\biggr] >0.
\end{equation}
The above inequality is not true for $s$ close to $0$ and larger values
of $k$. However, it is true
for $s>1$ and $k\in\{\lfloor s+1\rfloor, \lceil s+1\rceil\}$ as we
will now show.
Indeed, define, for $x\in[0,1]$,
\[
u_k(x) = \bigl[(1-x)2^{s/(s+1)}+x\bigr]^{(1+1/s)} -\biggl(2- \frac
{k}{k-1}x\biggr),
\]
and note that $u_k(j/k)$ is equal to the left-hand side of \eqref
{condPoiss}. Hence, if we show that when $s>1$ for both $k=\lfloor
s+1\rfloor$ and $k=\lceil s+1\rceil$, the function $u_k(x)> 0$ for all
$x\in(0,1)$, then we are done.
Differentiating $x\mapsto u_k(x)$ with respect to $x$ yields
\[
u^\prime_k(x)= -a[(1-x)2^{1/a}+x]^{a-1} (2^{1/a}-1)+ \frac{k}{k-1},
\]
where $a=(s+1)/s> 1$. The function $u_k^\prime$ is increasing as can
easily be seen from the second derivative
\[
u''_k(x)=a(a-1)(1-2^{1/a})^2[(1-x)2^{1/a}+x]^{a-2}>0.
\]
Hence, since $u_k(0)=0$, it suffices to show that $u^{\prime}_k(0)>0$,
for the two indicated values of~$k$.

\begin{claim*}\label{cl1} Fix $s>1$, then the statement $u^{\prime
}_k(0)>0$ is
true for both $k=\lfloor s+1\rfloor$ and $k=\lceil s+1\rceil$.
\end{claim*}

\begin{pf} Since $a=(s+1)/s$, the condition $u^{\prime
}_k(0)>0$ is equivalent to
\[
sk\bigl(2^{1/(s+1)}-1\bigr)>\bigl(k-(s+1)\bigr)\bigl(2-2^{1/(s+1)}\bigr).
\]
This inequality is trivially true for $k=\lfloor s+1\rfloor$, since
then the right-hand side is smaller than or equal to $0$, whereas the
left-hand side is positive.
We now turn to the case where $k=\lceil s+1\rceil$.
Since $2^{1/(s+1)}\ge1$ and $\lceil s+1\rceil-(s+1)\leq1$, the
right-hand side is bounded by 1, that is,
\[
\bigl(\lceil s+1\rceil-(s+1)\bigr)\bigl(2-2^{1/(s+1)}\bigr)\le1,\qquad s>1.
\]
A lower bound for the left-hand side on the domain $s>1$, is attained
in the limit as $s\downarrow1$ and equals $3(\sqrt{2}-1)=1.2426\ldots,$
that is,
\[
sk\bigl(2^{1/(s+1)}-1\bigr)=s\lceil s+1\rceil\bigl(2^{1/(s+1)}-1\bigr)\geq3\bigl(\sqrt
{2}-1\bigr),\qquad s>1.
\]
This shows that the above claim holds and hence that the Poisson
approximation holds both for $k=\lfloor s+1\rfloor$
and $k=\lceil s+1\rceil$.
\end{pf}
With the verification of the above claim the proof of Proposition~\ref{prop:poisson-dist-err} is complete.
\end{pf*}

\subsection{\texorpdfstring{The case $s\in\SSS$, the special set}{The case s in $\SSS$, the special set}}
\label{sec_bad}
In this section, we will prove Theorem~\ref{theo:main-s-in-ss}.
To this end, we fix $s_j\in\SSS$ and write $k=\lfloor s_j+1\rfloor$,
so that
$k+1=\lceil s_j+1\rceil$. Let $N_k^{(n)}=N_k^{(n)}(z_n(x))$
denote the number
of paths from 1 to $n$ of length $k$ and with weight at most $z_n(x)$,
with $z_n(\cdot)$ given by
\eqref{corr_zn}, and similarly we denote by $M_k^{(n)}=M_k^{
(n)}(z_n(y))$ the number of paths from 1 to $n$ of length $k+1$ and
with weight at most $z_n(y)$, where
$z_n(y)$ is given by the right-hand side of \eqref{corr_zn}, with $t$
replaced by $y$ and $k$ by $k+1$. Note that the change from $k$ to
$k+1$ is for many aspects irrelevant, because for $s=s_j$, we have
$g_s(k)=g_s(k+1)$. We are therefore, in particular, allowed to use the
same quantity $z_n(y)$ in the definition of $M_k^{(n)}$.
We show below that the total variation distance between
$N_k^{(n)}+M_k^{(n)}$ and a Poisson variable with mean
$\mu^{(n)}_k(x,y)=\expec[N_k^{(n)}+M_k^{(n)}]$ converges
to $0$ as $n\to\infty$; that is,
we show that
%
\begin{equation} \label{tot_var} d_{\mathrm{TV}}\bigl(N_k^{(n)}+M_k^{(n)},
\Poi\bigl(\mu^{ (n)}_k(x,y)\bigr)\bigr) \to0.
\end{equation}
Let us first prove that \eqref{tot_var} implies Theorem \ref
{theo:main-s-in-ss}.

\begin{pf*}{Proof of Theorem \protect\ref{theo:main-s-in-ss} assuming \protect\eqref{tot_var}}
The convergence in total variation in \eqref{tot_var} implies that
\begin{eqnarray*}
\prob\bigl(W_k(n)> z_n(x),W_{k+1}(n)> z_n(y)\bigr)
&=&\prob\bigl(N_k^{(n)}=0, M_k^{(n)}=0\bigr)
=\prob\bigl(N_k^{(n)}+M_k^{(n)}=0\bigr)\\
&\rightarrow& \prob\bigl(\Poi(\mu_k(x,y))=0\bigr),
\end{eqnarray*}
where, by \eqref{eqn:lambdak-asymp}, $\mu_k(x,y)=\lim_{n\rightarrow
\infty} \mu^{(n)}_k(x,y)=\lim_{n\to\infty} \expec[N_k^{
(n)}]+\expec[M_k^{(n)}]=\lambda_k(x)+\lambda_{k+1}(y)$, and where
we define
%
\begin{equation} \label{lambda-k-def} \lambda_l(z)=a_l \e^z, \qquad l\ge1,
z\in\Rbold.
\end{equation}
Thus, comparing with \eqref{duality},
%
\begin{eqnarray}\label{Poi-lim-ind}
&&\lim_{n\rightarrow\infty}
\prob\bigl(W_k(n)> z_n(x),W_{k+1}(n)> z_n(y)\bigr)
\nonumber
\\[-9pt]
\\[-9pt]
\nonumber
&&\quad = \lim_{n\rightarrow\infty
} \prob\bigl(W_k(n)> z_n(x)\bigr) \lim_{n\rightarrow\infty}\prob\bigl(W_{k+1}(n)>
z_n(y)\bigr),
\end{eqnarray}
and consequently, we see that the events
$\{W_k(n)> z_n(x)\}$ and $\{W_{k+1}(n)> z_n(y)\}$ are asymptotically
independent. It is then straightforward to conclude that
the minimum of the normalized pair $(W_k(n),W_{k+1}(n))$, where the
normalization is as in the left-hand side
of part (a) of Theorem~\ref{theo:main-s-in-ss} converges in
distribution to the minimum of the \textit{independent} pair
\[
\bigl(\Xi_k/(k-1),\Xi_{k+1}/k\bigr).
\]
The lower bound for $(\log n)^s W_k(n)$ of Theorem~\ref{lowerbound} and
the tightness
of $H_n$ (see (\ref{tightnessHn})) again completes the proof of part
(b), the hopcount part, and subsequently also part (a), of
Theorem~\ref{theo:main-s-in-ss}.
\end{pf*}

In order to prove \eqref{tot_var}, we again rely on the Poisson
approximation in \cite{BHJ92}. Set $\mathcal{T}_k(n)=\SSS_k(n)\cup\SSS
_{k+1}(n)$, the index set
of all paths from 1 to $n$ having either $k$ or $k+1$ edges, where, as
before, $k=\lfloor s_j+1\rfloor$. To denote that the length of a path
is equal to $k$, we give it a subscript
$k$ and write $\alpha_k$ for an element of $\SSS_k(n)$. For $\alpha
_k\in\SSS_k(n)$, we denote by
\[
I_{\alpha_k}=I_{\alpha_k}(z_n(x))= \mathbh{1}_{\{w_s(\alpha
_k)\leq z_n(x)\}},
\]
whereas for a path $\alpha_{k+1}\in\SSS_{k+1}(n)$, we define
\[
I_{\alpha_{k+1}}=I_{\alpha_{k+1}}(z_n(y))= \mathbh{1}_{\{
w_s(\alpha _{k+1})\leq z_n(y)\}},
\]
so that
\begin{eqnarray*}
p_{k}^{(n)}(x)&=&\prob\bigl(w_s(\alpha_k)\leq
z_n(x)\bigr)=F_k(z_n(x))\\[-2pt]
p_{k+1}^{(n)}(y)&=&\prob\bigl(w_s(\alpha_{k+1})\leq z_n(y)\bigr)=F_{k+1}(z_n(y)).
\end{eqnarray*}
Writing $\alpha$ for $\alpha_k$ or $\alpha_{k+1}$, we denote by
$\IIst(\alpha)\subseteq\mathcal{T}_k(n)$ the set of paths (not including~$\alpha$)
which have at least one edge in common to $\alpha$, and by
$\SSst(\alpha)\subseteq\mathcal{T}_k(n)$ the set of paths that do not
overlap on any edge with $\alpha$.
Finally, let
\[
Z_\alpha= \sum_{\beta_k\in\IIst(\alpha)}\mathbh{1}_{\{
w_s(\beta_k)\leq z_n(x)\}}
+\sum_{\beta_{k+1}\in\IIst(\alpha)}\mathbh{1}_{\{w_s(\beta
_{k+1})\leq z_n(y)\}}.
\]
The total variation distance in \eqref{tot_var} is bounded by
\begin{eqnarray*}
\label{eqn:dtv-bound2}
&&\frac{\sum_{\alpha\in\mathcal{S}_k(n)}[(p^{
(n)}_k(x))^2+p^{(n)}_k(x)\expec[Z_\alpha]+\expec[I_\alpha Z_\alpha]
]}{\mu^{(n)}_k(x,y)}
\\[-2pt]
&&\quad{}+
\frac{\sum_{\alpha\in\mathcal{S}_{k+1}(n)}[(p^{{
(n)}}_{k+1}(y))^2+p^{{(n)}}_{k+1}(y)\expec[Z_\alpha]+\expec
[I_\alpha Z_\alpha]
]}{\mu^{(n)}_k(x,y)}\nonumber.
\end{eqnarray*}
Since
\[
\mu^{(n)}_k(x,y)=
p_{k}^{(n)}(x)|\SSS_k(n)|+p_{k+1}^{(n)}(y)|\SSS_{k+1}(n)|
\ge\max\bigl\{p_{k}^{(n)}(x)|\SSS_k(n)|,p_{k+1}^{(n)}(y)|\SSS
_{k+1}(n)|\bigr\},
\]
we conclude from the proof of Proposition~\ref{prop:poisson-dist-err} that
\begin{eqnarray*}
\frac{\sum_{\alpha\in\mathcal{S}_k(n)}[(p^{
(n)}_k(x))^2+p^{(n)}_k(x)\expec[Z_\alpha]
]}{\mu^{(n)}_k(x,y)}
+
\frac{\sum_{\alpha\in\mathcal{S}_{k+1}(n)}[(p^{
(n)}_{k+1}(y))^2+p^{(n)}_{k+1}(y)\expec[Z_\alpha]
]}{\mu^{(n)}_k(x,y)}\to0.
\end{eqnarray*}
Hence, it remains to prove that
%
\begin{equation} \label{last-bound2} \frac{\sum_{\alpha\in\mathcal{S}_k(n)} \expec[I_\alpha Z_\alpha]+ \sum_{\alpha\in\mathcal{S}_{k+1}(n)} \expec[I_\alpha Z_\alpha] }{\mu^{(n)}_k(x,y)}\to0.
\end{equation}
We next decompose $\expec[I_\alpha Z_\alpha]$ into the part where
$\beta
$ has $k$ or $k+1$ edges, that is,
\[
\expec[I_\alpha Z_\alpha]=\sum_{\beta_k\in\IIst(\alpha)} \prob
(I_\alpha
=1,I_{\beta_k}=1)
+
\sum_{\beta_{k+1}\in\IIst(\alpha)} \prob(I_\alpha=1,I_{\beta_{k+1}}=1).
\]
By making this decomposition, as well as differentiating between
the number of edges of $\alpha$, the numerator in \eqref{last-bound2}
splits into 4 different double sums. The two double sums running over
the index sets $\alpha\in\mathcal{S}_k(n),\beta_k\in\IIst(\alpha)$
and $\alpha\in\mathcal{S}_{k+1}(n),\beta_{k+1}\in\IIst(\alpha)$ are
treated in the proof of Proposition~\ref{prop:poisson-dist-err}, apart from
the small change that $z_n(x)$ and $z_n(y)$ are now possibly
\textit{different}. Since $z_n(x)/z_n(y)\rightarrow1$, it is
straightforward to adapt the argument. Below, we will show that
%
\begin{equation} \label{doublesum} \frac{\sum_{\alpha\in\mathcal{S}_k(n)}\sum_{\beta_{k+1}\in\IIst(\alpha ) } \prob(I_\alpha
=1,I_{\beta_{k+1}}=1) }{\mu^{(n)}_k(x,y)}\to0.
\end{equation}
The terms with $\alpha\in\mathcal{S}_{k+1}(n)$ and $\beta_{k}\in\IIst
(\alpha)$
are identical, apart from the fact that $x$ and $y$ are interchanged.
Thus, \eqref{doublesum} completes the proof of \eqref{tot_var}.

To prove \eqref{doublesum}, we write, as in \eqref{eqn:st-sum-alpha},
\[
\sum_{\beta_{k+1}\in\IIst(\alpha) }
\prob(I_\alpha=1,I_{\beta_{k+1}}=1)
=\sum_{j=1}^{k-1} |\IIst_{k+1,j}(\alpha)|p_{k+1,j}^{(n)}(x,y),
\label{eqn:st-sum}
\]
where $\IIst_{k+1,j}(\alpha)\subset\mathcal{T}_k(n)$ consists of the
set of
paths of length $k+1$ which overlap with $\alpha$, which has length
$k$, in exactly $j$ edges, while
\[
p_{k+1,j}^{(n)}(x,y)=
\prob\bigl(X_{k,k}\leq z_n(x), X_{k+1,j} \leq z_n(y)\bigr),
\]
where, similarly as in the proof of Proposition~\ref{prop:poisson-dist-err},
we now write $X_{k,k}=\sum_{r=1}^k E_r^{-s},$ while
$X_{k+1,j}=\sum_{r=1}^j E_r^{-s}+ \sum_{r=j+1}^{k+1} \tilde
{E}_r^{-s}, 1\le j \le k-1$.

By adapting the Laplace method-type argument used in the proof of
Lemma~\ref{lemma:bi-bound}, it is readily verified that
for $z_1,z_2\downarrow0$ such that $\lim_{z_1\rightarrow0} z_2/z_1=1$,
and $k\geq3$ and $1\leq j\leq k-1$, we have
%
\begin{equation} \label{bik+1-bound} \prob(X_{k,k}\leq z_1,
X_{k+1,j}\leq z_2)=\exp\bigl(-z^{-p}_1[(k-j)\nu +j+1]^{p+1} \bigl(1+\mathrm{o}(1)\bigr)\bigr),
\end{equation}
where, as before, $\nu=2^{1/(p+1)}$.
By \eqref{zn-asympt} $(z_n(x))^{-p}\sim(k-1) \log n/(k^{1+p})$.
Further, since $s\in\mathcal{S}$, we have that $(z_n(y))^{-p}\sim\log
n/((k+1)^{1+p})
=(k-1) \log n/(k^{1+p}).$ Thus, $\lim_{n\rightarrow
\infty} z_n(y)/ z_n(x)=1,$ as required.

We conclude that in order for \eqref{doublesum} to hold, we need to
show that $n^{2k-j-1}p_{k+1,j}^{(n)}(t)\to0$,
or equivalently that
\[
\exp\biggl(-\log{n}\biggl[(k-1)
\biggl(\biggl(1-\frac{j}{k}\biggr)\nu+ \frac{j+1}{k}\biggr)^{p+1}-(2k-j-1)
\biggr]\biggr) \to0,
\qquad
1\le j \le k-1.
\]
This follows from the 
convexity of the function $x\mapsto x^{p+1}$ and the facts that $1<\nu
=2^{1/(p+1)}<2$ and $p>0$, since, for $1\le j \le k-1$,
\begin{eqnarray*}
\frac{(2k-j-1)}{k-1}&=&
\biggl(2-\frac{j-1}{k-1}\biggr)\leq\biggl(\biggl(1-\frac{j-1}{k}\biggr)\nu
+\frac{j-1}{k}\biggr)^{p+1}\\
&=&\biggl(\biggl(1-\frac{j}{k}\biggr)\nu+\frac{j+1}{k}+\frac{\nu-2}{k}\biggr)^{p+1}
< \biggl(\biggl(1-\frac{j}{k}\biggr)\nu+\frac{j+1}{k}\biggr)^{p+1}.
\end{eqnarray*}
This proves \eqref{doublesum}, and thus completes the proof of \eqref
{tot_var}.

\subsection{Multipoint distance limits}
\label{sec:multi-distance}
In this section, we indicate how to prove Corollary \ref
{corr:multi-point} for $2$ multipoint distances. The case for general
$m$ follows similarly.

More precisely, let $\tilde{W}_n^{(12)}, \tilde{W}_n^{(13)}$
denote the recentered and rescaled optimal
weights between~$1$ and $2$ and $1$ and $3$. Recall, for any fixed
$t\in\Rbold$, the function $z_n(t)$
from \eqref{corr_zn}, where we take $k=k^*(s)$. For $j=2,3$ and any
$t\in\Rbold$, let $N_k^{j,(n)}(z_n(t))$ denote the number of paths
between $1$ and $j$ having $k$ edges whose weight is less than $z_n(t)$.

The proof of Corollary~\ref{corr:multi-point} will be an adaptation of
the proof of Theorem~\ref{theo:main-s-in-ss} in Section~\ref{sec_bad}, and
we start by recalling some results we have proved and shall rely on.
Recall that we have already proved that, as $n\rightarrow\infty$,
\[
\lambda_k^{(n)}(t)=\expec\bigl[N_k^{j,(n)}(z_n(t))\bigr] \to\lambda_k(t),
\]
where $\lambda_k(t)=a_k \e^t$ is defined in \eqref{lambda-k-def} and
\[
d_{\mathrm{TV}}\bigl(N_k^{j,(n)}(z_n(t)), \Poi(\lambda_k(t))\bigr) \to0.
\]
For any fixed $x, y\in\Rbold$, define
%
\begin{equation} N_n^* = N_k^{2,(n)}(z_n(x))+ N_k^{3,(n)}(z_n(y)).
\end{equation}
Below, we shall show that
%
%
\begin{equation}
\label{eqn:sum-poi}
N_n^* \convd\Poi\bigl(\lambda_k(x)+\lambda_k(y)\bigr).
\end{equation}
Then the argument leading to \eqref{Poi-lim-ind} implies that
$\tilde{W}_n^{(12)}$ and $\tilde{W}_n^{(13)}$ are
asymptotically independent,
so that
\[
\lim_{n\to\infty}\prob\bigl(\tilde{W}_n^{(12)}>x, \tilde{W}_n^{
(13)}>y\bigr)\to\exp\bigl(-\lambda_k(x)-\lambda_k(y)\bigr),
\]
establishing the result we want. We next sketch how to prove \eqref
{eqn:sum-poi}.

\begin{pf*}{Sketch of proof of \protect\eqref{eqn:sum-poi}}
Fix any path $\alpha$ with $k$ edges between $1$ and $2$ and path
$\beta
$ with $k$ edges between $1$ and $3$. Since the argument is quite close
to the proof of Theorem~\ref{theo:main-s-in-ss}, we shall keep the discussion brief and focus on
the differences. We again rely on the total variation bound in Theorem~\ref
{theo:stein-chen-poi} that implies
\[
d_{\mathrm{TV}}\bigl(N_n^*, \Poi\bigl(\lambda_k^{(n)}(x)+\lambda_k^{(n)}(y)\bigr)\bigr)
\leq\frac{\mathrm{(I)}+\mathrm{(II)}+\mathrm{(III)} }{\lambda_k^{(n)}(x)+\lambda_k^{(n)}(y)}.
\]
Here,
\begin{eqnarray*}
\mathrm{(I)}&=&p_{k}^{(n)}(x)\lambda_k^{(n)}(x)
+p_{k}^{(n)}(y)\lambda_k^{(n)}(y),\\
\mathrm{(II)}&=& \lambda_k^{(n)}(x)\bigl(\expec\bigl[Z_\alpha^{(1,2)}\bigr]+\expec
\bigl[Z_\alpha^{(1,3)}\bigr]\bigr)
+\lambda_k^{(n)}(y)\bigl(\expec\bigl[Z_\beta^{(1,2)}\bigr]+\expec\bigl[Z_\beta
^{(1,3)}\bigr]\bigr),\\
\mathrm{(III)}&=&\expec\bigl[I_\alpha Z_\alpha^{(1,2)}+ I_\alpha Z_\alpha^{(1,3)}+
I_\alpha Z_\beta^{(1,2)}+ I_\alpha Z_\beta^{(1,3)}\bigr],
\end{eqnarray*}
where, as in \eqref{I-alpha-def},
\[
I_\alpha= \mathbh{1}_{\{|w_s(\alpha)|\leq z_n(x)\}},\qquad I_\beta=
\mathbh{1}_{\{|w_s(\beta)|\leq z_n(y)\}},
\]
while, writing $\IIst_{1,2}(\alpha)$ for the set of paths from $1$ to
$2$, which overlap with $\alpha$, and
$\IIst_{1,3}(\alpha)$ for the set of paths from $1$ to $3$, which
overlap with $\alpha$
(and similarly for $\beta$),
\[
Z_\alpha^{(1,2)} = \sum_{\gamma\in\IIst_{1,2}(\alpha)}\mathbh{1}_{\{ |w_s(\gamma)|\leq z_n(x)\}},
\quad
\mbox{and}
\quad
Z_\alpha^{(1,3)} = \sum_{\gamma\in\IIst_{1,3}(\alpha)}\mathbh{1}_{\{ |w_s(\gamma)|\leq z_n(y)\}},
\]
and similarly for $Z_\beta^{(1,2)}$ and $Z_\beta^{(1,3)}$.
Now, we have already shown that the terms (I) and (II) divided by
$\lambda_k^{(n)}(x)+\lambda_k^{(n)}(y)$ vanish as
$n\to\infty$. Thus, to complete the proof we just need to show that,
as $n\rightarrow\infty$,
\[
\frac{\expec[I_\alpha Z_\alpha^{(1,2)}]}{\lambda_k^{
(n)}(x)+\lambda_k^{(n)}(y)}\to0,
\]
as well as the corresponding other three terms of (III).
This is a minor adaptation of the proof of~\eqref{doublesum},
and we omit the details.
\end{pf*}

\section*{Acknowledgements}
The research of SB is supported by a UNC Research Council grant and
Junior Faculty Development award and NSF-DMS Grant 1105581.
SB would also like to thank the hospitality of \textsc{Eurandom}
where part of this work was done. The work of RvdH is supported
in part by Netherlands Organization for Scientific Research (NWO). We
thank Maren Eckhoff for
carefully proofreading the paper.
%

\printhistory

\end{document}